\documentclass[10pt]{article}
\usepackage{amssymb,amsmath,amsthm,latexsym}
\usepackage{amsfonts}
\usepackage[twlrm]{rawfonts}
\usepackage[english]{babel}
\usepackage[latin1]{inputenc}
\usepackage{fancyhdr}
\usepackage{indentfirst}
\usepackage{color}
\newtheorem{teo}{Theorem}[section]

\newtheorem{lem}{Lemma}[section]
\newtheorem{pro}{Proposition}[section]
\newtheorem{claim}{Claim}[section]
\newenvironment{dem}[1][Proof]{\noindent\textbf{#1.} }{\hfill \rule{0.5em}{0.5em}}

\newcommand{\N}{\mathbb{N}}

\newcommand{\R}{\mathbb{R}}

\begin{document}
	
\setlength{\baselineskip}{6.5mm} \setlength{\oddsidemargin}{8mm}
\setlength{\topmargin}{-3mm}

\title{\bf Existence of solution for a nonlocal problem in $\R^N$ via bifurcation theory}

\author{Claudianor O. Alves\thanks{C.O.Alves was partially supported by CNPq/Brazil 301807/2013-2 and INCT-MAT, coalves@dme.ufcg.edu.br}\, , Romildo N. de Lima\thanks{R.N. de Lima was partially supported by CAPES/Brazil, romildo@dme.ufcg.edu.br}\, and \,Marco A. S. Souto \thanks{M.A.S. Souto was partially supported by CNPq/Brazil
		305384/2014-7 and INCT-MAT, marco@dme.ufcg.edu.br}\,\,\,\,\,\,\vspace{2mm}
		\and {\small  Universidade Federal de Campina Grande} \\ {\small Unidade Acadêmica de Matemática} \\ {\small CEP: 58429-900, Campina Grande - PB, Brazil}\\}

\date{}

\maketitle

\begin{abstract}
	In this paper, we study the existence of solution for the following class of nonlocal problem,
$$
\left\{
\begin{array}{lcl}
-\Delta u=\left(\lambda f(x)-\int_{\R^N}K(x,y)|u(y)|^{\gamma}dy\right)u,\quad \mbox{in} \quad \R^{N}, \\
\displaystyle \lim_{|x| \to +\infty}u(x)=0,\quad u>0 \quad \text{in} \quad  \R^{N}, 
\end{array}
\right.
\eqno{(P)}
$$
where $N\geq3$, $\lambda >0, \gamma\in[1,2)$, $f:\R\rightarrow\R$ is a positive continuous function and $K:\R^N\times\R^N\rightarrow\R$ is a nonnegative function. The functions $f$ and $K$ satisfy some conditions, which permit to use Bifurcation Theory to prove the existence of solution for problem $(P)$\vspace{0.5cm}


\end{abstract}
	
\section{Introduction and main result}

The main goal of this paper is to study the existence of positive solutions for the following class of nonlocal problem 
$$
\left\{
\begin{array}{lcl}
-\Delta u=\left(\lambda f(x)-\int_{\R^N}K(x,y)|u(y)|^{\gamma}dy\right)u,\quad \mbox{in} \quad \R^{N}, \\
\displaystyle \lim_{|x|\to +\infty}u(x)=0,\quad u>0 \quad \text{in} \quad \R^{N}, 
\end{array}
\right.
\eqno{(P)}
$$
where $N\geq3$, $\lambda >0, \gamma\in[1,2)$, $f:\R^{N}\rightarrow\R$ is a positive continuous function and $K:\R^N\times\R^N\rightarrow\R$ is a nonnegative function. The functions $f$ and $K$ satisfy some technical conditions, which will be mentioned later on.

The motivation to study problem $(P)$ comes from the problem to model the behavior of a species inhabiting in a smooth bounded domain $\Omega\subset\R^N$, whose the classical logistic equation is given by 
\begin{equation} \label{1}
\left\{
\begin{array}{lcl}
-\Delta u=u(\lambda-b(x)u^{p}),\mbox{ in} \quad \Omega, \\
u=0, \quad  \mbox{on} \quad \partial\Omega,
\end{array}
\right.
\end{equation}
where $u(x)$ is the population density at location $x\in\Omega$, $\lambda\in\R$ is the growth rate of the species, and $b$ is a positive function denoting the carrying capacity, that is, $b(x)$ describes the limiting effect of crowding of the population.

Since (\ref{1}) is a local problem, the crowding effect of the population $u$ at $x$ only depends on the value of the population in the same point $x$. In \cite{Chipot}, for more realistic situations, Chipot has considered that the crowding effect depends also on the value of the population around of $x$, that is, the crowding effect depends on the value of integral involving the function $u$ in the ball $B_{r}(x)$ centered at $x$ of radius $r>0$. To be more precisely, in \cite{Chipot}, the following nonlocal problem has been studied   
\begin{equation} \label{2}
\left\{
\begin{array}{lcl}
-\Delta u=\left(\lambda-\int_{\Omega\cap B_{r}(x)}b(y)u^{p}(y)dy\right)u, \quad \mbox{in} \quad \Omega, \\
u=0, \quad \mbox{ on } \quad \partial\Omega,
\end{array}
\right.
\end{equation}
where $b$ is a nonnegative and nontrivial continuous function. After \cite{Chipot}, a special attention has been given for the problem 
\begin{equation} \label{3}
\left\{
\begin{array}{lcl}
-\Delta u=\left(\lambda-\int_{\Omega}K(x,y)u^{p}(y)dy\right)u,\mbox{ in } \quad \Omega, \\
u=0, \quad \mbox{ on } \quad  \partial\Omega,
\end{array}
\right.
\end{equation}
by supposing different conditions on $K$, see for example, Allegretto and Nistri \cite{Allegretto-Nistri}, Alves, Delgado, Souto and  Suárez \cite{Alves-Delgado-Souto-Suarez}, Chen and Shi \cite{Chen-Shi}, Corrêa, Delgado and Suárez \cite{Correa-Delgado-Suarez}, Coville \cite{Coville}, Leman, Méléard and Mirrahimi \cite{Leman-Meleard-Mirrahimi}, and Sun, Shi and Wang \cite{Sun-Shi-Wang} and their references. 

In \cite{Alves-Delgado-Souto-Suarez}, Alves, Delgado, Souto and  Suárez have considered the existence and nonexistence of solution for Problem (\ref{3}). In that paper, the authors have introduced a class $\mathcal{K}$, which is formed by functions $K:\Omega\times\Omega\rightarrow\R$ verifying:

$(i)$ $K\in L^{\infty}(\Omega\times\Omega)$ and $K(x,y)\geq0$ for all $x,y\in\Omega$.

$(ii)$ If $w$ is mensurable and $\int_{\Omega\times\Omega}K(x,y)|w(y)|^{p}|w(x)|^{2}dxdy=0$, then $w=0$ a.e. in $\Omega$.

\noindent Using Bifurcation Theory, and supposing that $K$ belongs to class $\mathcal{K}$, the following result has been proved 

\begin{teo} \label{T1}
	The problem $(3)$ has a positive solution if, and only if, $\lambda>\lambda_{1}$, where $\lambda_{1}$ is the first eigenvalue of problem
	\begin{equation*}
	\left\{
	\begin{array}{lcl}
	-\Delta u=\lambda u, \quad \mbox{in} \quad \Omega, \\
	u=0, \quad \mbox{on} \quad \partial \Omega.
	\end{array}
	\right.
	\end{equation*}
\end{teo} 

Motivated by \cite{Alves-Delgado-Souto-Suarez}, at least from a mathematical point of view, it seems to be interesting to ask if 
$(P)$ has a solution. Here, as in \cite{Alves-Delgado-Souto-Suarez}, we intend to use Bifurcation Theory, however we should be careful, because in the above paper, $\Omega$ is a smooth bounded domain, then it is possible to use compact embeddings for Sobolev spaces, and also, for Schauder spaces, like $H^{1}_{0}(\Omega)$ and $C^{2}(\overline{\Omega})$. Since  we are working in whole $\mathbb{R}^{N}$, we need to show new estimates, and to this end, our inspiration comes from of some articles due to Edelson and Rumbos \cite{Edelson-Rumbos1, Edelson-Rumbos2}, where the Bifurcation Theory has been used to study the existence of solution for a problem of the type
$$
\left\{
\begin{array}{lcl}
-\Delta u+u=\lambda f(x)(u+h(u)),\quad \mbox{in} \quad \R^{N}, \\
\displaystyle \lim_{|x|\to +\infty}u(x)=0,\quad u>0 \quad \text{in} \quad \R^{N}, 
\end{array}
\right.
\eqno{(Q)}
$$ 
where $f$ and $h$ are continuous function verify some technical conditions. Here we would like to point out that a lot of estimates in our paper are totally different from those used in \cite{Edelson-Rumbos1, Edelson-Rumbos2}, because in the present work the problem is nonlocal, while there the problem is local.

In the present article, we assume that $K:\R^N\times\R^N\rightarrow\R$ is a continuous functions verifying the following conditions:

\noindent $(K_{0})$\, There is $P\in C^{+}_{rad}(\R^N,\R) \cap L^{1}(\R^N) $ such that 
$$
0\leq K(x,y)\leq f(x)P(y)^{\gamma/2}Q(x,y), \quad \forall (x,y) \in \R^N \times \R^N,
$$ 
where, $C^{+}_{rad}(\R^N,\R)=\left\{g \in C(\R^N,\R)\,:\,g \quad \mbox{is positive and radially symmetric} \right\}$ and $Q:\R^N\times\R^N\rightarrow\R$ is a mensurable function verifying:
\begin{description}
	\item[$(Q_{1})$] $M=\displaystyle \sup_{x \in \mathbb{R}^{N}}|Q(x,.)|_{\frac{2}{2-\gamma}}<+\infty$.
	\item[$(Q_{2})$] Given $\varepsilon>0$, there exist $R,L>0$ such that
	$$\int_{|y|\leq L}Q(x,y)^{2/2-\gamma}dy<\varepsilon,\quad \forall x\in B^{c}_{R}(0).$$
\end{description}

\noindent $(K_{1})$ If $w$ is mensurable and $ \int_{\R^N\times\R^N}K(x,y)|w(y)|^{\gamma}|w(x)|^{2}dxdy=0$, then $w=0$ a.e. in $\R^N$.

\noindent From $(K_0)-(K_1)$, it follows that  $K\in L^{1}(\R^N\times\R^N)$. 

Related to function $f$, we assume that

\noindent $(f_{1})$ $f:\R\rightarrow\R$ is continuous such that $0<f(x)\leq P(x), \quad \forall x\in\R^N.$

\noindent $(f_{2})$ There exists $q>N/2$, such that
$$
\sup_{x\in\R^N}|f|_{L^{q}(B_{2}(x))}<+\infty.
$$

By taking $P\in C^{+}_{rad}(\R^N,\R) \cap L^{1}(\Omega)$ and $Q(x,y)=g(y-x)$, where $g\in L^{2/2-\gamma}(\R^N)$ the function   
$$
K(x,y)=f(x)P(y)^{\gamma/2}Q(x,y)
$$ 
verifies the conditions $(K_{0})$ and $(K_{1})$. 

Our main result is the following 

\begin{teo} \label{T2}
Assume that $(K_0)-(K_1)$ and $(f_1)-(f_2)$ hold. Then, problem $(P)$ has a positive solution if, and only if, $\lambda>\lambda_{1}$, where $\lambda_{1}$  is the first eigenvalue of the linear problem
$$
	\left\{
	\begin{array}{l}
	-\Delta u=\lambda f(x)u, \quad \mbox{in} \quad \R^N,  \\
	\displaystyle \lim_{|x|\to +\infty}u(x)=0.
	\end{array}
	\right.
	\eqno{(AP)}
$$
\end{teo}

The paper is organized as follows: In Section 2, we have shown some properties of the nonlocal term. In Section 3, we have defined two compact operators which are crucial in our approach, while in Section 4, we prove the Theorem \ref{T2}.

\noindent \textbf{Notations}
\begin{itemize}
	
	\item $\omega_{N}$ is the volume of the unit ball in $\R^N$.
	
	\item $\Gamma$ is the fundamental solution of Laplace equation in $\R^N$.
	
	\item $\chi_{B}$ is the characteristic function of $B$.
	
	\item $B_{r}(x)$ denotes the ball centered at the $x$ with radius $r>0$ in $\R^N$.
	
	\item $L^{s}(\R^N)$, for $1\leq s\leq\infty$, denotes the Lebesgue space with usual norm denoted by $|u|_s$. 
	
	\item $L^{2}_{P}(\R^N)$ denotes the class of real valued Lebesgue measurable functions $u$ satisfying
	$$
	\int_{\R^N}P(x)|u(x)|^{2}dx<\infty,
	$$
	which is a Hilbert space endowed with the inner product
	$$
	(u,v)_{2,P}=\int_{\R^N}P(x)u(x)v(x)dx, \quad \forall  u,v\in L^{2}_{P}(\R^N).
	$$
	The norm associated with this inner product will denote by $|\,\,\,\,\,|_{2,P}$.
	
	\item $D^{1,2}(\R^N)$ denotes the Sobolev space endowed with inner product
	$$
	(u,v)_{1,2}=\int_{\R^N}\nabla u\nabla vdx, \forall u,v\in D^{1,2}(\R^N).
	$$
The norm associated with this inner product will denote by $\|\,\,\,\,\,\|_{1,2}$. In \cite{Edelson-Rumbos1}, it was proved that the embedding  $ D^{1,2}(\R^N) \hookrightarrow L^{2}_{P}(\R^N)$ is compact.

\item We denote by $E$, the Banach space given by
$$
E=\left\{u\in C(\R^N);\lim_{|x|\rightarrow\infty}u(x)=0\right\},
$$
endowed with the norm $|\,\,\,\,\,|_\infty$.  A simple computation gives that the embedding  $ E\hookrightarrow L^{2}_{P}(\R^N)$ is continuous. 	
	\item If $u$ is a mensurable function, we denote by $u^{+}$ and $u^{-}$  the positive and negative part of $u$ respectively, which are given by
	$$
	u^{+}=\max\{u,0\} \quad \mbox{and} \quad u^{-}=\max\{-u,0\}.
	$$
	
\end{itemize}

\section{The nonlocal term}

In the sequel, we will show some properties of the operator $\phi:L^{2}_{P}(\R^N) \to L^{1}(\R^N)$, given by \linebreak $\phi(u):=\phi_{u}$, where
$$
\phi_{u}(x):=\int_{\R^N}K(x,y)|u(y)|^{\gamma}dy,
$$
which is well defined, because we are assuming $(K_0)-(K_1)$. Using the definition of $\phi$, we see that $u \in D^{1,2}(\R^N)$ is a solution for $(P)$ if, and only if, it is a positive solution of 
$$
-\Delta u+\phi_{u}u=\lambda f(x)u,\quad\mbox{in}\quad\R^{N}. \eqno{(EP)}
$$

In the sequel, we show some important properties of the operator $\phi$ for future reference

\begin{lem} \label{phi}
The operator $\phi$ satisfies the following properties: \\
\noindent ${ (\phi_{1})}$ $\phi_{tu}=t^{\gamma}\phi_{u}$, \quad $\forall (u,t) \in E \times [0,+\infty)$;

\noindent ${(\phi_{2})}$ \, 
$$
|\phi_{u}(x)|\leq MP(x)|u|^{\gamma}_{2,P} \quad \forall (u,x) \in  L^{2}_{P}(\R^N) \times \R^N \quad \mbox{and} \quad |\phi_{u}(x)|\leq MP(x)|u|^{\gamma}_{\infty} \quad \forall (u,x) \in E \times \R^N;
$$ 
\noindent $(\phi_3)$ 
For each $u \in  L^{2}_{P}(\R^N)$,
$$
\lim_{|x|\to +\infty} \frac {\phi_u(x)}{f(x)}=0 .
$$
\noindent ${(\phi_{4})}$ 
$$
|\phi_{u}|_{1}\leq M|P|_{1}|u|^{\gamma}_{2,P} \quad \forall u \in  L^{2}_{P}(\R^N) \quad \mbox{and} \quad |\phi_{u}|_{1}\leq M|P|_{1}|u|^{\gamma}_{\infty} \quad \forall u \in E.  
$$
\noindent ${(\phi_{5})}$ $\phi:E\longrightarrow L^{1}(\R^N)$ is continuous, that is, 
$$
u_n \to u \quad \mbox{in} \quad E \Longrightarrow \phi_{u_n} \to \phi_u \quad \mbox{in} \quad L^{1}(\R^{N}).
$$
\noindent $(\phi_6)$ \, Let $(u_n) \subset L^{2}_{P}(\R^N) $ be a sequence and $u \in E$ such that $u_n(x) \to u(x)$ a.e in $\R^N$. Then
$$
\liminf_{n \to +\infty}\phi_{u_n}(x) \geq \phi_{u}(x), \quad \forall x \in \R^N.
$$ 
\noindent $(\phi_7)$ \, Let $(u_n) \subset E$ be a sequence and $u \in E$ such that $u_n \to u$ in $E$. Then, given $\epsilon >0$, there is $n_0 \in \mathbb{N}$ such that
$$
|\phi_{u_n}(x)-\phi_u(x)|\leq \epsilon P(x), \quad \forall n \geq n_0 \quad \mbox{and} \quad \forall x \in \mathbb{R}^{N}.
$$
\noindent $(\phi_8)$ \, If $u \in C^{1}(\R^N) \cap D^{1,2}(\R^N)$ is a  nontrivial solution of $(EP)$ and  $u \geq 0$ ( resp. $u \leq 0$ ), then $u>0$ \linebreak ( resp. $u<0$).
\end{lem}

\begin{dem} \\
\noindent { $(\phi_1)$:} \, This property is an immediate consequence of the definition of $\phi_u$.\\
\noindent { $(\phi_2)$:} \, From $(K_0)$ and $(f_1)$, for any $u \in L^{2}_{P}(\R^N)$,
$$
|\phi_{u}(x)|\leq P(x)\int_{\R^N}P(y)^{\frac{\gamma}{2}}Q(x,y)|u(y)|^{\gamma}dy.
$$
Then, using H\"older inequality with exponents $p=\frac{2}{\gamma}$ and $p'=\frac{2}{2-\gamma}$, we get
$$
|\phi_{u}(x)|\leq P(x)|Q(x,y)|_{\frac{2}{2-\gamma}}|u|_{2,P}^{\gamma} \leq M P(x)|u|_{2,P}^{\gamma} \quad \forall u \in  L^{2}_{P}(\R^N),
$$
where $M=\displaystyle \sup_{x \in \mathbb{R}^{N}}|Q(x,.)|_{\frac{2}{2-\gamma}}$ was fixed in $(Q_1)$. The last inequality combined with  the continuous embedding $E \hookrightarrow  L^{2}_{P}(\R^N) $ gives 
$$
|\phi_{u}(x)|\leq M P(x)|u|_{\infty}^{\gamma} \quad \forall u \in  E \quad \mbox{and} \quad \forall x \in \R^N.
$$
\noindent {  $(\phi_3)$:} \, Repeating the same arguments explored in $(\phi_2)$, for each $L>0$ we have the inequality below
$$
|\phi_{u}(x)|\leq f(x)\left[ \left(\int_{|y|\leq L}Q(x,y)^{\frac{2}{2-\gamma}}dy \right)^{\frac{2-\gamma}{2}}|u|_{2,P}^{\gamma} + M\left(\int_{|y|>L }P(y)|u(y)|^{\gamma}dy \right)^{\frac{\gamma}{2}} \right], \quad \forall x \in \R^N.
$$
Combining the fact that $u \in L^{2}_{P}(\R^N)$ with $(Q_2)$, given $\epsilon >0$, there are $R,L$ such that 
$$
|\phi_{u}(x)|\leq \epsilon f(x) \quad \mbox{for} \quad |x| \geq R,
$$
showing $(\phi_3)$. \\
\noindent {  $(\phi_4)$:}\, This property follows from $(\phi_2)$, because $P \in L^{1}(\R^N)$.\\
\noindent {  $(\phi_5)$:}\, An immediate consequence from Lebesgue's Theorem together with $(\phi_2)$.\\ 
\noindent {  $(\phi_6)$:}\, As $K(x,y)$ is nonnegative, the property is obtained applying Fatous' Lemma. \\
\noindent {  $(\phi_7)$:}\, Using $(f_1)$ and the definitions of $\phi_{u_n}$ and $\phi_u$, we get
$$
|\phi_{u_n}(x)-\phi_u(x)|\leq P(x)\int_{\R^N}P(y)^{\frac{\gamma}{2}}Q(x,y)||u_n(y)|^{\gamma}-|u(y)|^{\gamma}|dy,
$$
and so,
$$
|\phi_{u_n}(x)-\phi_u(x)|\leq MP(x)|P|_{1}^{\frac{\gamma}{2}}|||u_n|^{\gamma}-|u|^{\gamma}|_\infty, \quad \forall x \in \R^N \quad \mbox{and} \quad \forall n \in \mathbb{N}.  
$$
As $u_n \to u$ in $E$, we have that
$$
||u_n|^{\gamma}-|u|^{\gamma}|_\infty \to 0 \quad \mbox{as} \quad n \to +\infty,
$$
showing that given $\epsilon$, there is $n_0 \in \mathbb{N}$ such that
$$
|\phi_{u_n}(x)-\phi_u(x)|\leq \epsilon P(x), \quad \forall x \in \R^N \quad \mbox{and} \quad n \geq n_0.
$$
\noindent {  $(\phi_8)$:}\, Immediate consequence of the maximum principles.
\end{dem}

\subsection{The weight $P$}

The aim of this section is to study the existence and regularity of some linear problems, which will be used in the proof of some lemmas later on. 

To begin with, if $P \in C^{+}_{rad}(\R^N,\R)$ and $\varphi$ is the weak solution of the problem 
$$
	\left\{
	\begin{array}{l}
	-\Delta u=P(x), \quad \mbox{in} \quad \R^N,  \\
	\displaystyle \lim_{|x|\to +\infty}u(x)=0,
	\end{array}
	\right.
	\eqno{(LP)}
$$
then $\varphi\in D^{1,2}_{rad}(\R^N)\cap C^{2}(\R^N)$ and 
\begin{equation} \label{E1}
-(r^{N-1}\varphi'(r))'=r^{N-1}P(r),\mbox{ for $r>0$ and }\lim_{r\rightarrow\infty}\varphi(r)=0.
\end{equation}

The lemma below shows the  behavior of $\varphi$ at infinite. A similar result has been proved \cite{Brezis-Kamin}, but with a different argument. 
\begin{lem} \label{primeirophi}
	The function $\varphi$ is decreasing, positive and
	$$
	\lim_{r \to +\infty}r^{N-2}\varphi(r)=\frac{|P|_{1}}{\omega_{N}(N-2)}.
	$$
\end{lem}

\begin{dem}
	Indeed, by (\ref{E1}), 
	$$
	-r^{N-1}\varphi'(r)=\int_{0}^{r}s^{N-1}P(s)ds, \quad \forall r>0.
	$$
	As $P$ is positive,  it follows that $\varphi'(r)<0$ for $r>0$, showing that $\varphi$ is decreasing. Moreover, 
	$$
	\lim_{r\longrightarrow\infty}-r^{N-1}\varphi'(r)=\int_{0}^{\infty}s^{N-1}P(s)ds=\frac{|P|_{1}}{\omega_{N}}.
	$$
	We also have 
	$$
	0<\varphi(r)=-\int_{r}^{\infty}\varphi'(s)ds\leq\frac{|P|_{1}}{\omega_{N}}\int_{r}^{\infty}s^{1-N}ds=\frac{|P|_{1}}{\omega_{N}(N-2)}r^{2-N},
	$$
	loading to
	\begin{equation} \label{E2}
	\limsup_{r \to +\infty}r^{N-2}\varphi(r)\leq \frac{|P|_{1}}{\omega_{N}(N-2)}.
	\end{equation}
	On the other hand, given $\varepsilon>0$, there exists $r_{0}>0$ such that
	$$
	-r^{N-1}\varphi'(r)\geq\frac{|P|_{1}-\varepsilon}{\omega_{N}}, \quad \mbox{for} \quad r>r_{0}.
	$$
	Therefore, if $r>r_{0}$,
	$$
	\varphi(r)=-\int_{r}^{\infty}\varphi'(s)ds\geq\frac{|P|_{1}-\varepsilon}{\omega_{N}}r^{2-N}.
	$$
Since $\varepsilon$ is arbitrary, we ensure that
	\begin{equation} \label{E3}
	\liminf_{r \to +\infty}r^{N-2}\varphi(r)\geq \frac{|P|_{1}}{\omega_{N}(N-2)}.
	\end{equation}
Now, the lemma follows combining (\ref{E2}) and (\ref{E3}). 
\end{dem}

The last lemma combined with some arguments found \cite[pages 225-226]{Edelson-Rumbos1} permits to conclude that 
\begin{equation} \label{varphi}
\varphi(x)=-\int_{\R^N}\Gamma(x-y)P(y)dy, \quad \forall x \in \R^N.
\end{equation}

\subsection{Some regularity results}

Let $F:\R^{N}\rightarrow\R$ be a continuous function verifying the following condition:
$$
|F(x)|\leq c_{0}P(x), \quad \forall x \in \R^N, \eqno{(H)}
$$
for some positive constant $c_0$. 

Since $w(x)\equiv 1\in L^{2}_{P}(\R^N)$, we can guarantee that the functional $\Psi:D^{1,2}(\R^N) \to \R$ given by
$$
\Psi(v):=\int_{\R^N}F(x)vdx
$$
is continuous, because  
$$
\Psi(v)=\int_{\R^N}F(x)vdx\leq c_{0}\int_{\R^N}P(x)vwdx\leq c_{0}|v|_{2,P}|P|_{1}^{1/2}\leq C|P|_{1}^{1/2}\|v\|_{1,2}.
$$

Using the Riesz's Theorem, there exists unique  $u\in D^{1,2}(\R^N)$ such that
$$
\int_{\R^N}\nabla u\nabla vdx=\int_{\R^N}F(x)vdx, \quad \forall v\in D^{1,2}(\R^N) \quad \mbox{and} \quad \|u\|_{1,2}\leq C|P|_{1}^{1/2}.
$$
Furthermore, by regularity theory , $u\in D^{1,2}(\R^N)\cap C^{1}(\R^N)$, and it is a weak solution of 
$$
-\Delta u=F(x) \quad \mbox{in} \quad \R^N.
$$ 

Using the above notation, we are able to prove the following result

\begin{pro} \label{F} 
Assume that $F$ satisfies the condition $(H)$. Then, there exists unique $u\in C^{1}(\R^N)\cap D^{1,2}(\R^N)$ with
	$$
	\int_{\R^N}\nabla u\nabla vdx=\int_{\R^N}F(x)vdx,\quad \forall v\in D^{1,2}(\R^N)
	$$
	and
	$$
	|u(x)|\leq\frac{c_{0}|P|_{1}}{\omega_{N}(N-2)}|x|^{2-N},\quad \forall x\in\R^N.
	$$
	Moreover, for each $R>0$,
	$$
	\|\nabla u\|_{C(B_R)}\leq\left[\frac{N}{R}\|u\|_{C(B_{2R})}+\frac{12R}{N-2}\|F\|_{C(B_{2R})}\right].
	$$
\end{pro}
\begin{dem}
In the sequel, we denote by $\rho:\R^N \to \R$ the function given by
$$
\rho(x)=
\left\{
\begin{array}{l}
e^{\frac{1}{|x|^{2}-1}}, \quad \mbox{if} \quad |x| <1, \\
0, \quad \mbox{if} \quad |x| \geq 1.
\end{array}
\right.
$$
It is well known that $\rho \in C^{\infty}_{0}(\R^N)$ with $supp \rho \subset \overline{B}_1(0)$. Using the function $\rho$, for each 
for each $n \in \mathbb{N}$, we set 
$$
u_{n}(x):=\int_{\R^N}\rho_{n}(x-y)u(y)dy,
$$
where $\rho_n(x)=Cn^{N}\rho(n x)$ with $C=(\int_{\R^N}\rho(y)\,dy)^{-1}$. Applying some results found \cite{Brezis}, we know that  
$(u_n)$ and $(\frac{\partial u_{n}}{\partial x_i})$ converge uniformly in compact of $\R^N$ for $u$ and $\frac{\partial u}{\partial x_i}$ respectively, for all $i \in \{1,...,N\}$. Moreover, fixing
$$
F_{n}(x)=\int_{\R^N}\rho_{n}(x-y)F(y)dy,
$$
we derive that $u_n$ verifies the equality below in the classical sense  
$$
-\Delta u_{n}=F_{n}(x), \quad \mbox{in} \quad \R^N.
$$
By priori estimates, see \cite{Gilbarg-Trudinger}, for each $R>0$:
$$
\|\nabla u_{n}\|_{C(B_{R})}\leq\left[\frac{N}{R}\|u_{n}\|_{C(B_{2R})}+\frac{12R}{N-2}\|F_{n}\|_{C(B_{2R})}\right], \quad \forall n \in \mathbb{N}.
$$
Passing to the limit $n \to +\infty$, and recalling that $(F_n)$ converges uniformly for $F$ in compacts, we deduce that 
$$
\|\nabla u\|_{C(B_{R})}\leq\left[\frac{N}{R}\|u\|_{C(B_{2R})}+\frac{12R}{N-2}\|F\|_{C(B_{2R})}\right].
$$
Using the Green's function $G_{R}$ of the ball $B_{R}$, with $R>|x|$,
$$
u_{n}(x)=-\int_{B_{R}}G_{R}(x,y)F_{n}(y)dy+\frac{(R^2-|x|^2)}{N\omega_{N}R}\int_{\partial B_{R}}\frac{u_{n}(\xi)}{|x-\xi|^{N}}d\sigma_{\xi}.
$$
Passing to the limit $n \to +\infty$, it follows that
$$
u(x)=-\int_{B_{R}}G_{R}(x,y)F(y)dy+\frac{(R^2-|x|^2)}{N\omega_{N}R}\int_{\partial B_{R}}\frac{u(\xi)}{|x-\xi|^{N}}d\sigma_{\xi},\mbox{ for} \quad R>|x|.
$$
Proceeding as in \cite[pages 225-226]{Edelson-Rumbos1}, we get
$$
u(x)=-\int_{\R^{N}}\Gamma(x-y)F(y)dy.
$$
Gathering the last equality with (\ref{varphi}), we find
$$
|u(x)|\leq \int_{\R^N}|\Gamma(x-y)| |F(y)|dy\leq c_{0}\int_{\R^N}|\Gamma(x-y)|P(y)dy=c_{0}|\varphi(x)|\leq \frac{c_{0}|P|_{1}}{\omega_{N}(N-2)}|x|^{2-N}. 
$$
The proof of the lemma is complete. \end{dem}

\section{A linear solution operator}

In this section, we study the existence and properties of an important operator, which will use to prove the existence of solution for problem $(P)$.

In what follows, we fix $f\in C(\R^N)$ with $0<f(x)\leq P(x)$. Then, for each $v\in L^{2}_{P}(\R^N)$, the compact embedding $ D^{1,2}(\R^N)\hookrightarrow L^{2}_{P}(\R^N)$ together with Riesz's Theorem yields there is a unique solution $u \in D^{1,2}(\R^N)$ of the problem
$$
	\left\{
	\begin{array}{l}
	-\Delta u=f(x)v, \quad \mbox{in} \quad \R^N,  \\
	\displaystyle \lim_{|x|\to +\infty}u(x)=0.
	\end{array}
	\right.
	\eqno{(WLP)_v}
$$

From this, we can define a {\it solution operator} $S:L^{2}_{P}(\R^N)\rightarrow L^{2}_{P}(\R^N)$ such that $S(v)=u$, where $u$ is the unique solution of the above weight linear problem $(WLP)_v$. By using well known arguments, $S$ is a compact self-adjoint operator, then by spectral theory, there exists a complete orthonormal basis $\{u_{n}\}$ of $L^{2}_{P}(\R^N)$, and a corresponding sequence of positive real numbers $\{\lambda_{n}\}$, with $\lambda_{n}\rightarrow\infty$ as $n\rightarrow\infty$, such that
$$
0<\lambda_{1}\leq\lambda_{2}\leq...\leq \lambda_n \leq .....
$$
and
$$
-\Delta u_{n}=\lambda_{n}f(x)u_{n},\mbox{ in $\R^N$}.
$$
Moreover, using Lagrange multiplier, it is possible to prove the following characterization for $\lambda_1$ 
$$
\lambda_{1}=\inf_{v\in D^{1,2}(\R^N) \setminus \{0\}}\frac{\int_{\R^N}|\nabla v|^{2}dx}{\int_{\R^N}f(x)|v(x)|^{2}dx}.
$$
The above identity is crucial to show that $\lambda_{1}$ is a simple eigenvalue, and that a corresponding eigenfunction $\varphi_{1}$ can be chosen positive in $\R^N$. Moreover, we also have the following property
\begin{equation} \label{comportamento u1}
\liminf_{|x| \to \infty}|x|^{N-2}\varphi_{1}(x)>0.
\end{equation}
The above limit is a consequence of the lemma below

\begin{lem} \label{comportamento u11} Let $u \in E \cap D^{1,2}(\R^N)$ be a positive function and $R>0$ such that  
$$
\int_{\R^N}\nabla u \nabla \psi \, dx \geq 0 \quad \forall \psi \in D_0^{1,2}(\R^N), \quad supp \ \psi \subset B^{c}_R(0) \quad \mbox{and} \quad \psi \geq 0.
$$
Then,
$$
\liminf_{|x| \to +\infty}|x|^{N-2}u(x)>0.
$$
\end{lem}
\begin{dem} First of all, as $u$ is a positive continuous function, there is $\varepsilon >0$ such that
$$
u(x)\geq\varepsilon>0,\quad \forall x\in \overline{B}_{R}(0).
$$
Setting $w(x)=u(x)-\frac{\varepsilon R^{N-2}}{2}|x|^{2-N}$ and 
$$
\tilde{w}(x)=\left\{
\begin{array}{l}
0, \quad \mbox{if} \quad |x| \leq R, \\
w^{-}(x), \quad \mbox{if} \quad |x| >R, 
\end{array}
\right.
$$
we have that  $\tilde{w} \in D_0^{1,2}(\R^N), supp \ \tilde{w} \subset B^{c}_R(0)$ and $\tilde{w} \geq 0$. Then, we must have
$$
\int_{\R^N}\nabla w\nabla \tilde{w} \ dx=\int_{\R^N}\nabla u \nabla \tilde{w} dx\geq0,
$$
implying that
$$
-\int_{\R^N \setminus B_R(0)}|\nabla w^{-}|^{2}dx\geq0.
$$
Then, $w^{-}=0$  in $\R^N\setminus B_{R}(0)$, from where it follows that 
$$
u(x)\geq\frac{\varepsilon R^{N-2}}{2}|x|^{2-N} \quad \mbox{for} \quad |x|\geq R, 
$$
finishing the proof of the lemma. 
\end{dem}

Since $E \subset L^{2}_{P}(\R^N)$, we intend to prove that $S:E \to E$ is well defined and it is a linear compact operator. 
To see why, we will consider the subspace $E_0$ of $E$ given by
$$
E_0:=\left\{u\in C^{1}(\R^N);\quad \sup_{x\in\R^N}[|x|^{N-2}u(x)]<\infty\right\},
$$
endowed with the following norm
$$
\|u\|:=\sup\left\{|x|^{N-2}|u(x)|:x\in\R^N\right\}.
$$

Using the space $E_0$, we claim that $S(E) \subset E_0$. Indeed, for each $v\in E$, considering  $F(x)=f(x)v(x)$, the Proposition \ref{F} ensures the existence of a unique $u\in C^{1}(\R^N)\cap D^{1,2}(\R^N)$ verifying
\begin{equation}
\int_{\R^N}\nabla u\nabla wdx=\int_{\R^N}f(x)vwdx, \quad \forall w\in D^{1,2}(\R^N)
\end{equation}
and
\begin{equation} \label{EstCima}
|u(x)|\leq \frac{|P|_{1}}{\omega_{N}(N-2)}|v|_{\infty}|x|^{2-N},\quad \forall x\in\R^N.
\end{equation}
Therefore, $S(v)=u \in E_0$. As $v \in E$ is arbitrary, we can guarantee that $S(E) \subset E_0$.

Next, we show an important relation involving $E_0$ and $E$, which will  use to show $S:E \to E$ is compact.
\begin{lem} \label{convergenciaK}
	Let $(u_n)$ be a bounded sequence in $E_0$. If for each compact $A \subset \R^N$, the sequence $(\|u_n\|_{C^{1}(A)})$ is also bounded, then $(u_n)$ admits a convergent subsequence in $E$.
\end{lem}

\begin{dem}
	Note that, from boundedness of $(u_n)$ in $E_0$, there exists $R_{1}>0$ such that
	$$
	|u_{n}(x)-u_{m}(x)|\leq M|x|^{2-N}<1,\quad \mbox{for} \quad |x|\geq R_{1} \quad \mbox{and} \quad \forall n \in \mathbb{N}.
	$$
	On the other hand, using the hypothesis that $(u_n)$ is bounded $C^{1}(\overline{B}_{R_1}(0))$, it follows that  
	$$
	|u_{n}(x)-u_{n}(y)|\leq M|x-y|,\quad \forall x,y\in\overline{B}_{R_{1}} \quad \mbox{and} \quad \forall n\in\N.
	$$
Applying the Arzelá-Ascoli's Theorem, there exists $\N_{1}\subset\N$, such that $(u_{n})_{n\in\N_{1}}$ is a Cauchy's sequence in $C(\overline{B}_{R_{1}})$ with 
$$
|u_{n}(x)-u_{m}(x)|<1, \quad \mbox{for} \quad |x|> R_{1} \quad \mbox{and} \quad n,m\in \N_{1}.
$$
Repeating the above arguments, there exists  $R_{2}>R_{1}$, such that
$$
|u_{n}(x)-u_{m}(x)|\leq M|x|^{2-N}<1/2, \quad \mbox{for} \quad |x|> R_{2} \quad \forall n,m\in\N_{1}.
$$
Once $(u_n)$ is bounded $C^{1}(\overline{B}_{R_2}(0))$, we derive that 
$$
|u_{n}(x)-u_{n}(y)|\leq M|x-y|, \quad \forall x,y\in\overline{B}_{R_{2}} \quad \mbox{and} \quad n\in\N_{1}.
$$
Applying again Arzelá-Ascoli's Theorem, there exists $\N_{2}\subset\N_{1}$, such that $(u_{n})_{n\in\N_{2}}$ is a Cauchy's sequence in $C(\overline{B}_{R_{2}}(0))$ and 
$$
|u_{n}(x)-u_{m}(x)|<1/2, \quad \mbox{for} \quad |x|> R_{2} \quad \mbox{and} \quad n,m\in \N_{2}.
$$
Repeating the above argument, we will find an increasing sequence $(R_k) \subset \R$ with $R_k \to +\infty$, and sets $\N\supseteq\N_{1}\supseteq\N_{2}\supseteq...\supseteq\N_{k}\supseteq...$ such that   
$$
|u_{n}(x)-u_{m}(x)|<1/k,\mbox{ always that $|x|> R_{k}$, with $n,m\in\N_{k}$},
$$
Thereby, there is a subsequence of $(u_n)$, still denoted by itself, such that, given $\varepsilon>0$, there exist $R>0$ and $n_{0}\in\N$ verifying 
$$
|u_{n}(x)-u_{m}(x)|<\varepsilon,\quad \mbox{for} \quad |x|>R \quad \mbox{and} \quad n,m\geq n_{0}.
$$
On the other hand, the boundedness of $(u_n)$ in $C^{1}(\overline{B}_{R_2}(0))$ permits to assume, changing the subsequence if necessary, the inequality below 
$$
|u_{n}(x)-u_{m}(x)|<\varepsilon \quad \mbox{for} \quad n,m\geq n_{0} \quad \mbox{and} \quad x\in \overline{B}_{R}(0).
$$ 
Therefore,
$$
|u_{n}-u_{m}|_{\infty}<\varepsilon,\mbox{ for all $n,m\geq n_{0}$}.
$$
Thus, $(u_{n})$ is a Cauchy's subsequence in $E$. Once $E$ is a Banach space, $(u_{n})$ is convergent in $E$, finishing the proof of the lemma.
\end{dem}

Now, we are ready to prove the compactness of $S$ 
\begin{lem} \label{compacidadeS}
	The operator $S:E\rightarrow E$ is compact.
\end{lem}
\begin{dem}
	Let $(v_n)$ be a bounded sequence in $E$ and $u_n=S(v_n)$. By (\ref{EstCima}), $u_n \in E_0$ and
$$
\|u_n\|\leq\frac{|v_n|_{\infty}|P|_{1}}{\omega_{N}(N-2)}, \quad \forall n \in \N.
$$
Moreover, considering $F_n(x)=f(x)v_n(x)$ and fixing $R>0$, the Proposition \ref{F} guarantees that 
	$$
	\|\nabla u_n\|_{C(B_R)}\leq\left[\frac{N}{R}\|u_n\|_{C(B_{2R})}+\frac{12R}{N-2}\|F_n\|_{C(B_{2R})}\right], \quad \forall n \in \N.
	$$
Using $(f_2)$ and bootstrap argument, we know that $(u_n)$ is also bounded in $C(\overline{B}_{2R}(0))$. As $(F_n)$ is also bounded in $C(\overline{B}_{2R}(0))$, because $(v_n)$ is bounded in $E$,  the right side of the last inequality is bounded. Thereby, we can apply the Lemma \ref{convergenciaK} to infer that $(u_n)$ possesses a convergent subsequence in $E$, and the lemmas follows. 
\end{dem}

Next, we show a result, which will be used in the proof of Theorem \ref{T2}.
\begin{lem} \label{positividade}
Let  $u\in E$ be a positive solution of $u=\lambda_{1}S(u)$ and $\sigma,R>0$ satisfying 
	$$
	|x|^{N-2}u(x)\geq\sigma,\quad \mbox{for} \quad |x|\geq R.
	$$
Let $v\in E$ be a weak solution of
$$
-\Delta v+b(x)v=\lambda f(x)v,\quad \mbox{in} \quad \R^N,
$$
where $b$ is a continuous funtions. Then, there exists $\varepsilon>0$ such that: if $|\lambda-\lambda_1|+|u-v|_{\infty}<\varepsilon$ and  $|b(x)|\leq \varepsilon P(x)$, for all $x\in\R^N$, the function $v$ is positive and $2|x|^{N-2}v(x)\geq\sigma$ for $|x|\geq R$.
\end{lem}
\begin{dem}
Indeed, for $w=u-v$ we have that
$$
\int_{\R^N}\nabla w\nabla\psi dx=\int_{\R^N}F(x)\psi,\quad \forall \psi\in D^{1,2}(\R^N),
$$
where 
$$
F(x)=(\lambda_{1}-\lambda)f(x)v+\lambda_{1}f(x)w-b(x)w+b(x)u, \quad \forall x \in \R^N.
$$ 
Using the hypotheses, it follows that 
$$
|F(x)|\leq C\varepsilon P(x), \quad \forall x \in \R^N.
$$ 
Thus, by Proposition \ref{F}, choosing $\varepsilon>0$ small enough, we see that 
$$
|x|^{N-2}|w(x)|\leq \frac{C\varepsilon |P|_{1}}{\omega_{N}(N-2)}<\frac{\sigma}{2} \quad \forall x \in \R^N,
$$
and so,
$$
2|x|^{N-2}|w(x)|\leq \sigma, \quad \forall x\in\R^N. 
$$
Consequently
$$
|x|^{N-2}v(x)\geq|x|^{N-2}u(x)-|x|^{N-2}w(x)\geq\frac{\sigma}{2},\quad \mbox{for} \quad |x|\geq R,
$$
showing that $v$ is positive for $|x|\geq R$. Now, for  $|x|\leq R$, decreasing $\varepsilon$ if necessary, the positiveness of $u$ gives that $v$ is also positive for $|x| \leq R$, finishing the proof of the lemma. 
\end{dem}

\subsection{ A nonlinear compact operator} 
In this subsection, we will study the properties of another compact operator, which will appear in our study. 

For each $v\in E$, using the notations of Section 2, there is $C>0$ such that  
$$
|-\phi_{v}(x)v(x)|\leq C|v|_{\infty}^{\gamma+1}P(x), \quad \forall x \in \R^N.
$$ 
Thus, applying Proposition \ref{F} with $F(x)=-\phi_{v}(x)v(x)$, there exists unique $u\in E_0 \cap D^{1,2}(\R^N)$ verifying 
\begin{equation} \label{phiu}
\int_{\R^N}\nabla u\nabla wdx=\int_{\R^N}-\phi_{v}(x)vwdx,\quad \forall w\in D^{1,2}(\R^N).
\end{equation}
Moreover, we also have 
$$
|u(x)|\leq\int_{\R^N}|\Gamma(x-y)| |F(y)|dy\leq C|v|_{\infty}^{\gamma+1}\int_{\R^N}\Gamma(x-y) P(y)dy\leq C|v|_{\infty}^{\gamma+1}|\varphi(x)|,
$$
where $\varphi$ was given in (\ref{varphi}). Once $\varphi$ is bounded, we get
\begin{equation} \label{E4}
|u(x)|\leq C|v|_{\infty}^{\gamma+1}, \quad \forall x \in \mathbb{R}^{N}.
\end{equation}

From the previous analysis, we can define the nonlinear operator $G:E\rightarrow E_0\subset E$, given by $G(v)=u$, where $u \in E_0 \cap D^{1,2}(\R^N)$ is the unique solution of (\ref{phiu}).  

The lemma below establishes that $G$ is a compact operator. Since the proof this fact follows with the same type of arguments explored in the proof of Lemma  \ref{compacidadeS}, we omit its proof.

\begin{lem} \label{compacidadeG}
	The operator $G:E\rightarrow E$ is compact.
\end{lem}
Using the definition of $G$, (\ref{E4}) yields  
$$
|G(v)|_{\infty}\leq C|v|_{\infty}^{\gamma+1}, \quad \forall v \in E,
$$
from where it follows that 
\begin{equation} \label{G1}
\lim_{|v|_{\infty}\to +\infty}\frac{G(v)}{|v|_{\infty}}=0,
\end{equation}
that is,
\begin{equation} \label{G2}
G(v)=o(|v|_{\infty}).
\end{equation}

\section{Proof of Theorem 2}

Using the definitions of $S$ and $G$, it easy to check that $(\lambda,u) \in \R \times D^{1,2}(\R^N)$ solves $(P)$ if, and only if, 
\begin{equation} \label{equacao}
u=F(\lambda,u):=\lambda S(u)+G(u).
\end{equation}
In the sequel, we will apply the following result due to Rabinowitz \cite{Rabinowitz}

\begin{teo}\textbf{(Global bifurcation)}
	Let $E$ be a Banach space. Suppose that $L$ is a compact linear operator
	and $\lambda^{-1}\in \sigma(L)$ has odd algebraic multiplicity . If $\Psi$ is a compact operator satisfying
	$$
	\lim_{\|u\|\to 0}\frac{\Psi(u)}{\|u\|}=0,
	$$
	then  the set
	$$
	\Sigma=\overline{\{(\lambda,u)\in\R\times E:u=\lambda L(u)+\Psi(u),u\neq0\}}
	$$
	has a closed connected component $\mathcal{C}=\mathcal{C}_{\lambda}$ such that $(\lambda,0)\in\mathcal{C}$ and
	
	(i) $\mathcal{C}$ is unbounded in $\R\times E$, or
	
	(ii) there exists $\hat{\lambda}\neq\lambda$, such that $(\hat{\lambda},0)\in\mathcal{C}$ and $\hat{\lambda}^{-1}\in\sigma(L)$.
\end{teo}

In what follows, we will apply the above theorem with $L=S$ and $\Psi=G$. By the previous results, we know that there is a first positive eigenfunction $\varphi_{1}$ associated with $\lambda_{1}$. Moreover, $\lambda_{1}^{-1}$ is an eigenvalue of $S$ with  multiplicity equal to 1. From global bifurcation theorem, there exists a closed connected component $\mathcal{C}=\mathcal{C}_{\lambda_{1}}$ of solutions for $(P)$, which satisfies $(i)$ or $(ii)$. We claim that $(ii)$ does not occur. In order to show this claim, we need of the lemma below

\begin{lem} \label{sinal}
	There exists $\delta>0$ such that, if $(\lambda,u)\in\mathcal{C}$ with $|\lambda-\lambda_1|+|u|_{\infty}<\delta$ and $u\neq0$, then $u$ has defined signal, that is,
$$
u(x)>0, \quad \forall x\in\R^N \quad \mbox{or} \quad u(x)<0,\quad \forall x\in\R^N.
$$
\end{lem}

\begin{dem}
It is enough to prove that for any two sequences $(u_n) \subset E$ and $\lambda_{n}\to \lambda_1$ satisfying 
	$$
	u_n\neq0,\quad |u_n|_{\infty}\to 0 \quad \mbox{and } \quad u_{n}=F(\lambda_n,u_n)=\lambda_{n}S(u_n)+G(u_n),
$$
$u_n$ has defined signal for $n$ large enough. 

Setting $w_n=u_n/|u_n|_{\infty}$, we have that 
	\begin{equation} \label{E6}
	w_n=\lambda_{n}S(w_n)+\frac{G(u_n)}{|u_n|_{\infty}}=\lambda_{n}S(w_n)+o_{n}(1).
	\end{equation}
From compactness of the operator $S$, we can assume that $(S(w_n))$ is convergent. Then, $w_n\to w$ in $E$, for some $w \in E$ with $|w|_\infty =1$. Consequently, 
$$
w=\lambda_{1}S(w)
$$
or equivalently,
$$
-\Delta w=\lambda_{1} f(x)w,\quad \mbox{ in } \quad \R^{N}.
$$
Thereby, $w$ is an eigenfunction associated with $\lambda_1$. Then, 
$$
w(x)>0 \quad \forall x\in\R^N \quad \mbox{ or } \quad w(x)<0 \quad \forall x\in\R^N.
$$
Without loss of generality, we assume that $w(x)>0$ for all $x\in\R^N$. 	As $w$ is the limit of $(w_{n})$ in $E$, by Lemma \ref{positividade}, it follows that 
$$
w_{n}(x)>0 \quad \forall x\in\R^N,
$$
for $n$ large enough. Once  $u_n$ and $w_n$ has the same signal, we have that $u_n$ is also positive, which is the desired conclusion. 
\end{dem}

It is easy to check that if $(\lambda,u)\in\Sigma$, the pair $(\lambda,-u)$ is also in $\Sigma$. In the lemma below, we show that $\mathcal{C}$ can be decomposed  into two important sets. 

\begin{lem} Consider the sets 
$$
\mathcal{C}^{+}=\{(\lambda,u)\in\mathcal{C}:u(x)>0,\quad \forall x\in\R^N \}\cup\{(\lambda_{1},0)\}
$$
and
$$
\mathcal{C}^{-}=\{(\lambda,u)\in\mathcal{C}:u(x)<0,\quad \forall x\in\R^N\}\cup\{(\lambda_{1},0)\}.
$$
Then, 
\begin{equation} \label{Componente}
\mathcal{C}=\mathcal{C}^{+}\cup\mathcal{C}^{-}.
\end{equation}
Moreover, note that   
$\mathcal{C}^{-}=\{(\lambda,u)\in\mathcal{C}:(\lambda,-u)\in\mathcal{C}^{+}\}$, $\mathcal{C}^{+}\cap\mathcal{C}^{-}=\{(\lambda_{1},0)\}$ and $\mathcal{C}^{+}$ is unbounded if, and only if, $\mathcal{C}^{-}$ is also unbounded.
\end{lem}
\begin{dem} \, In what follows, we fix
$$
\mathcal{C}^{\pm}=\{(\lambda,u)\in\mathcal{C}:u^{\pm}\not=0 \}.
$$
To prove (\ref{Componente}), it is enough to show that $\overline{\mathcal{C}^{\pm}} = \emptyset$. Arguing by contradiction, if $\overline{\mathcal{C}^{\pm}} \not= \emptyset$, as  $\mathcal{C}$ is a connected set in $\R \times E$, we must have 
$$
\left( \mathcal{C}^{+} \cup \mathcal{C}^{-} \right) \cap \overline{\mathcal{C}^{\pm}} \not= \emptyset.
$$ 
Therefore, there is a solution $(\lambda, u)$ of $(P)$  and sequences $(\lambda_n, u_n) \subset \mathcal{C}^{+} \cup \mathcal{C}^{-} $ and $(s_n,w_n) \subset \mathcal{C}^{\pm}$ such that 
$$
\lambda_n, s_n \to \lambda \quad \mbox{in} \quad \mathbb{R} , \quad u_n \to u \quad \mbox{in} \quad E \quad \mbox{and} \quad w_n \to u \quad \mbox{in} \quad E. 
$$ 
Consequently $u \geq 0$ in $\mathbb{R}^{N}$ or $u \leq 0$ in $\mathbb{R}^{N}$, then by Lemma \ref{sinal}, $u \not =0$. Supposing that $u \geq 0$ and $u \not= 0$, the property $(\phi_8)$ ensures that $u(x)>0$ in $\mathbb{R}^{N}$.  Now, by $(\phi_3)$, there is $R>0$ such that
$$
-\Delta u=(\lambda f(x)-\phi_u(x))u\geq 0, \quad  \mbox{for} \quad |x| \geq R,
$$
in the weak sense, that is, 
$$
\int_{\R^N}\nabla u \nabla \psi \, dx \geq 0 \quad \forall \psi \in D^{1,2}(\R^N), \quad supp \ \psi \subset B^{c}_R(0) \quad \mbox{and} \quad \psi \geq 0.
$$
Applying Lemma \ref{comportamento u11}, we get
\begin{equation} \label{P2}
\liminf_{|x| \to +\infty}|x|^{2-N}u(x)>0.
\end{equation}
Now, setting
$$
F_n(x)=(\lambda_{n}-\lambda)f(x)w_n+\lambda f(x)(w_n-u)+(\phi_{w_n}-\phi_u)w_n+\phi_u (w_n-u), 
$$
given $\epsilon >0$, the proprieties $(\phi_2)$ and $(\phi_7)$ guarantee that there is $n_0 \in \mathbb{N}$ such that
$$
|F_n(x)|< \epsilon P(x), \quad \forall n \geq n_0 \quad \mbox{and} \quad \forall x \in \mathbb{R}^{N}.
$$
The above inequality combined with (\ref{P2}) permits to repeat the same arguments used in the proof of Lemma \ref{positividade} to conclude that $w_n$ is positive for $n$ large enough, obtaining a contradiction. Thereby,  $\overline{\mathcal{C}^{\pm}} = \emptyset$, finishing the proof. 
\end{dem}

Now, we are able to prove that $(ii)$ does not hold.  
\begin{lem}
	$\mathcal{C}^{+}$ is unbounded.
\end{lem}
\begin{dem}
	Suppose by contradiction that $\mathcal{C}^{+}$ is bounded. Then, $\mathcal{C}$ is also bounded. From global bifurcation theorem, there exists $(\hat{\lambda},0)\in\mathcal{C}$ , where $\hat{\lambda}\neq\lambda_{1}$ and $\hat{\lambda}^{-1}\in\sigma(S)$.
Hence, there exists $(u_n)$ in $E$ and $\lambda_{n}\rightarrow\hat{\lambda}$, such that
$$
u_{n}\neq0,\quad |u_{n}|_{\infty}\rightarrow0\quad\mbox{and}\quad u_{n}=F(\lambda_{n},u_{n}).
$$
Considering $w_{n}=u_{n}/|u_{n}|_{\infty}$, we know that (\ref{E6}) is also satisfied. Moreover, as in the proof of Lemma \ref{sinal}, passing to a subsequence if necessary, $(w_n)$ converges to $w$ in $E$, which is a nontrivial solution of the problem 
$$
-\Delta w=\hat{\lambda}f(x)w,\quad \mbox{in} \quad \R^N,
$$
showing that $w$ is an eigenfunction related to $\hat{\lambda}$. Since that $\hat{\lambda}\neq\lambda_{1}$, $w$ must change signal. Then, for $n$ large enough, $w_n$ must change signal, implying that $u_n$ also changes signal, which is an absurd, 	because $(\lambda_{n},u_{n})\in\mathcal{C}^{+}$ or $(\lambda_{n},u_{n})\in\mathcal{C}^{-}$.
\end{dem}

From previous lemma, the connected component $\mathcal{C}^{+}$ is unbounded. Now, our goal is to show that this component
intersects any hiperplane $\{\lambda\}\times E$, for $\lambda>\lambda_{1}$.  To see this, we need of the following priori estimate

\begin{lem}\textbf{(A priori estimate)} \label{priori} For any $\Lambda>0$, there exists $r>0$ such that, if $(\lambda,u)\in\mathcal{C}^{+}$ and $\lambda \in [0,\Lambda]$, then $|u|_{\infty}\leq r$. 
\end{lem}

\begin{dem}
	We start the proof with the following claim:
	\begin{claim}
		For any $\Lambda>0$, there exists $r>0$ such that, if $(\lambda,u)\in\mathcal{C}^{+}$ and $\lambda \in [0,\Lambda]$, we must have $\|u\|_{1,2}\leq r$. Consequently, by Sobolev embedding, there is $r_1>0$ such that $|u|_{2^{*}} \leq r_1$.
	\end{claim}
\noindent Indeed, if the claim does not hold, there are $(u_{n})$ in $D^{1,2}(\R^N)$ and $(\lambda_{n})\subset[0,\Lambda]$ such that $\|u_{n}\|_{1,2}\to \infty$ and $u_{n}=F(\lambda_{n},u_{n})$. Considering $w_{n}=u_{n}/\|u_{n}\|_{1,2}$, it follows that,
	\begin{equation} \label{E7}
	\int_{\R^N}\nabla w_{n}\nabla\psi dx+\int_{\R^N}\phi_{u_{n}}(x)w_{n}\psi dx=\lambda_{n}\int_{\R^N}f(x)w_{n}\psi dx,\quad \forall \psi\in D^{1,2}(\R^N).
		\end{equation}
	Once $(w_{n})$ is bounded in $D^{1,2}(\R^N)$, without loss of generality, we can suppose that there is $w\in D^{1,2}(\R^N)$, such that $w_{n}\rightharpoonup w$ in $D^{1,2}(\R^N)$. Consequently, for some subsequence, the Sobolev embedding and $(\phi_5)$ combine to give
\begin{equation} \label{E8}
w_{n}(x) \to w(x)\quad \mbox{a.e. in} \quad \R^N \quad \mbox{and} \quad \liminf_{n \to \infty}\phi_{w_{n}}(x) \geq \phi_{w}(x) \quad \forall x \in \R^N.
\end{equation}
Setting $\psi=u_{n}/\|u_{n}\|^{\gamma+1}_{1,2}$ as a function test in (\ref{E7}) and using $(\phi_1)$, we get 
$$
\frac{1}{\|u_{n}\|_{1,2}^{\gamma}}+\int_{\R^N}\phi_{w_{n}}(x)w_{n}^{2}dx=\frac{\lambda_{n}}{\|u_{n}\|_{1,2}^{\gamma}}\int_{\R^N}f(x)w_{n}^{2}dx,\quad \forall n\in\N.
$$
Hence,  
$$
\lim_{n \to \infty}\int_{\R^N}\phi_{w_{n}}(x)w_{n}^{2}dx=0.
$$
Then, the Fatous' Lemma together with (\ref{E8}) loads to   
$$
0\leq\int_{\R^N}\phi_{w}(x)w^{2}(x)dx \leq \lim_{n}\int_{\R^N}\phi_{w_{n}}(x)w_{n}^{2}dx=0,
$$
that is,
$$
\int_{\R^N\times\R^N}K(x,y)|w(y)|^{\gamma}|w(x)|^{2}dydx=0.
$$
Therefore, from $(K_1)$, it follows that $w\equiv0$, and so, $w_{n} \to 0$ in $L^{2}_{P}(\R^N)$.
Now, fixing $\psi=w_{n}$ as a test function in (\ref{E7}), we obtain
$$
\int_{\R^N}|\nabla w_{n}|^{2}dx+\int_{\R^N}\phi_{u_{n}}(x)w_{n}^{2}dx=\lambda_{n}\int_{\R^N}f(x)w_{n}^{2}dx,
$$
from where it follows that
$$
\int_{\R^N}|\nabla w_{n}|^{2}dx\leq \Lambda\int_{\R^N}P(x)w_{n}^{2}dx\to 0.
$$
Thus,
$$
\|w_{n}\|^{2}_{1,2} \to 0
$$
which is an absurd, because $\|w_{n}\|_{1,2}=1$, for all $n\in\N$, which proves the claim. 	

In order to obtain a priori estimate, we need of a good estimate from above for the norm $|\,\,\,\,|_{\infty}$. To this end, we will use the lemma below, whose proof will be omitted, because it is a small modification of the iteration Moser, similar to that found in \cite{Chabrowski-Szulkin}. 

\begin{lem} \label{Moser}
	Let $h:\R^N\longrightarrow\R$ be a nonnegative mensurable function verifying 
	$$
	\sup_{x\in\R^N}|h|_{L^{q}(B_{2}(x))}<\infty
	$$
	with $q>N/2$ and $v\in D^{1,2}(\R^N)$ be a weak solution of the problem
	\begin{equation}
	-\Delta v+b(x)v=H(x,v),\quad \mbox{in} \quad \R^N,
	\end{equation}
	where $b:\R^N \to \R$ is a nonnegative function and $ H:\R^N\times\R\longrightarrow\R$ is a continuous function verifying,
	$$
	|H(x,s)|\leq h(x)|s|,\quad \forall (x,s) \in \R^N \times \R.
	$$
	Then, $v \in E$ and there exists a constant $C:=C(q,h)>0$ such that
	$$
	|v|_{\infty}\leq C|v|_{2^{*}}.
	$$
\end{lem}	
To conclude the proof of Lemma \ref{priori}, it is sufficient to apply the Lemma \ref{Moser} with
$$
b(x)=\phi_{u_n}(x), \quad H(x,s)=\lambda_n f(x)s \quad \mbox{and} \quad h(x)=\Lambda f(x).
$$ 
\end{dem}

\noindent {\bf Conclusion of the proof of Theorem \ref{T2}}

To finalize the proof of Theorem \ref{T2}, we must show that there is no solution for $(P)$ when $\lambda < \lambda_{1}$. Indeed, arguing by contradiction, if $(\lambda,u)$ is a solution of $(P)$, taking $\psi=\varphi_{1}$ as test function in $(P)$, we derive that 
$$
\int_{\R^N}\nabla u\nabla\varphi_{1}dx+\int_{\R^N}\phi_{u}u\varphi_{1}dx=\lambda\int_{\R^N}f(x)u\varphi_{1}dx
$$
from where it follows that,
$$
\lambda_{1}\int_{\R^N}f(x)u\varphi_{1}<\lambda\int_{\R^N}f(x)u\varphi_{1},
$$
showing that  $\lambda_{1}<\lambda$, and the proof is complete.

\end{document}